\documentclass[a4,11pt]{amsart}%
\usepackage{amsmath,amssymb,amsfonts,enumerate,amsthm, amscd,}
\usepackage[all]{xy}
\usepackage{graphicx}
\usepackage{amsmath}
\usepackage{amsfonts}
\usepackage{amssymb}
\usepackage{hyperref}%

\newcommand{\N}{\mathbb{N}}

\newtheorem{thm}{Theorem}[section]
\newtheorem{cor}[thm]{Corollary}
\newtheorem{lem}[thm]{Lemma}

\newtheorem{defn}[thm]{Definition}

\newtheorem{exam}[thm]{Example}
\newtheorem{rem}[thm]{Remark}

\def\proof{{\parindent0pt {\bf Proof.\ }}}
\newcommand{\field}[1]{\mathbb{#1}}

\newcommand{\Z }{\field{Z}}

\newcommand{\p}{\mathfrak{p}}

\theoremstyle{definition}

\theoremstyle{remark}

\begin{document}
	
	\bibliographystyle{amsplain}

	\title[$S$-Prime and $S$-maximal ideals in trivial ring extensions]{$S$-Prime and $S$-maximal ideals in trivial ring extensions of commutative rings}
	\author{Hwankoo Kim}
	\address{Hwankoo Kim \\Division of Computer Engineering, Hoseo University, Asan 31499, Republic of Korea.}
	\email{hkkim@hoseo.edu}

	\author{Najib Mahdou}
	\address{Najib Mahdou\\Department of Mathematics, Faculty of Science and Technology of Fez, Box 2202, University S.M. Ben Abdellah Fez, Morocco.}
	\email{mahdou@hotmail.com}
	
	\author{El Houssaine Oubouhou}
	\address{El Houssaine Oubouhou, Department of Mathematics, Faculty of Science and Technology of Fez, Box 2202, University S.M. Ben Abdellah Fez, Morocco.}
	\email{hossineoubouhou@gmail.com}


	\begin{abstract}
This paper explores the study of $S$-prime and $S$-maximal ideals in the context of trivial ring extensions $A \ltimes M$. Through counterexamples, we demonstrate that $S$-prime (resp., $S$-maximal) ideals in $A \ltimes M$ are not necessarily homogeneous, and a homogeneous $S$-prime (resp., $S$-maximal) ideal does not necessarily have the form $P \ltimes M$, where $P$ is an $S_0$-prime (resp., $S_0$-maximal) ideal of $A$.

Moreover, we characterize the conditions under which an ideal $J$ (not necessarily homogeneous) in the trivial ring extension $A \ltimes M$ is $S$-prime (resp., $S$-maximal). Additionally, we demonstrate that all $S$-prime (and consequently $S$-maximal) ideals in $A \ltimes M$ are of the form $P \ltimes M$, where $P$ is an $S_0$-prime ideal of $A$, if and only if $M$ is an $S_0$-divisible $A$-module.

As an application, we explore the transfer of the concepts of compactly $S$-packed rings, coprimely $S$-packed rings and   $S$-$pm$-rings to the trivial ring extension.

These results provide significant insights into the relation between $S$-primality and $S$-maximality in trivial ring extensions, contributing to a deeper understanding of ideal theory in this context. This work not only enriches the theoretical framework of ring structures but also advances the broader field of algebraic theory through practical examples and applications.

	\end{abstract}
	
	\maketitle

	
	

	
\section{Introduction}

Throughout this article, we consider only commutative rings with identity. The symbol $R$ will always denote such a ring. The notion of prime ideals and their generalizations play a significant role in Commutative Algebra and Algebraic Geometry. They are used to characterize a wide variety of rings and have applications in other areas such as General Topology, Graph Theory, Cryptology, and more. A proper ideal $\mathfrak{p}$ of $R$ is said to be a {\em prime ideal} if, whenever $a b \in \mathfrak{p}$ for some $a, b \in R$, then $a \in \mathfrak{p}$ or $ b \in \mathfrak{p}$. The set of all prime ideals and maximal ideals of $R$ are denoted by $\operatorname{Spec}(R)$ and $\operatorname{Max}(R)$, respectively. 

One of the most important methods for generalizing classical rings and modules involves the use of multiplicative subsets $S$ of $R$. A nonempty subset $S$ of $R$ is called a {\em multiplicative subset} of $R$ if it satisfies the following conditions: $1 \in S$, $0 \notin S$, and for every $a, b \in S$, we have $ab \in S$. Notable examples of multiplicative subsets include the set of regular elements $\operatorname{reg}(R)$ and $S_P = R \setminus P$ for each $P \in \operatorname{Spec}(R)$ (see \cite{Gilmer}). Given a multiplicative subset $S$ of $R$, its saturation is defined as the set $S^* = \{x \in R \mid xy \in S \text{ for some } y \in R\}$, which is itself a multiplicative subset of $R$ satisfying $S \subseteq S^*$. Various authors have studied this concept extensively (see, for instance, \cite{ABT, H1, H2, KKL, MMZ, Sartinian}).

At the beginning of this century, in 2002, Anderson and Dumitrescu introduced the notion of an {\em $S$-Noetherian ring} in \cite{AD}. A ring $R$ is said to be $S$-Noetherian if, for every ideal $I$ of $R$, there exists a finitely generated sub-ideal $K$ of $I$ such that $sI \subseteq K$ for some $s \in S$. Several fundamental results, including Cohen’s Theorem, the Eakin-Nagata Theorem, and the Hilbert Basis Theorem for $S$-Noetherian rings, were established in \cite{AD}. However, the choice of $s \in S$ ensuring $sI \subseteq K$ in the definition of $S$-Noetherian rings is not uniform. Consequently, Qi et al. \cite{QKWCZ} introduced the notion of {\em uniform $S$-Noetherian rings} and obtained results such as the Eakin-Nagata-Formanek Theorem and the Cartan-Eilenberg-Bass Theorem for these rings.

Recently, Sevim et al. \cite{S} studied the concept of {\em $S$-prime ideals}, a generalization of prime ideals, and used it to characterize integral domains, specific prime ideals, fields, and $S$-Noetherian rings. An ideal $\mathfrak{p}$ of $R$ is called an {\em $S$-prime ideal} if $\mathfrak{p} \cap S = \emptyset$ and there exists a fixed $s \in S$ such that  if, whenever $a b \in \mathfrak{p}$ for some $a, b \in R$, then $s a \in \mathfrak{p}$ or $s b \in \mathfrak{p}$ \cite{HM, S}. The set of all $S$-prime ideals of $R$ is denoted by $\operatorname{Spec}_S(R)$. Similarly, an ideal $\mathfrak{p}$ of $R$ is called an {\em $S$-maximal ideal} if $\mathfrak{p} \cap S = \emptyset$ and there exists a fixed $s \in S$ such that whenever $\mathfrak{p} \subseteq Q$ for some ideal $Q$ of $R$, then either $sQ \subseteq \mathfrak{p}$ or $Q \cap S \neq \emptyset$ \cite{szariski}. The set of all $S$-maximal ideals of $R$ is denoted by $\operatorname{Max}_S(R)$. More recently, Yildiz et al. constructed a topology on $\operatorname{Spec}_S(R)$, generalizing the classical Zariski topology \cite{szariski}, and investigated topological properties such as connectedness, compactness, and separation axioms.

Let $A$ be a ring and $M$ an $A$-module. The following ring construction, known as the {\em trivial extension} of $A$ by $M$ (also called the idealization of $M$), is denoted by $A \ltimes M$. It is the ring whose additive structure is that of the external direct sum $A \oplus M$, with multiplication defined by
\[
(r_1, m_1) (r_2, m_2) := (r_1r_2, r_1m_2 + r_2m_1)
\]
for all $r_1, r_2 \in A$ and all $m_1, m_2 \in M$. This construction is sometimes referred to using different terminology or notations, such as the idealization $A( + ) M$. While the exact origin of this construction is unclear, its first systematic study appears in Nagata’s book \cite{nagata}. The purpose of idealization is to embed $M$ into a commutative ring $R$ such that the structure of $M$ as an $A$-module is preserved as an ideal of $R$.

The trivial ring extension is useful for:
\begin{enumerate}
    \item Reducing results concerning submodules to the ideal case,
    \item Generalizing results from rings to modules, and
    \item Constructing examples of commutative rings with zero-divisors.
\end{enumerate}
The fundamental properties of trivial ring extensions are summarized in \cite{gl, Huck}. Primarily, trivial ring extensions have proven valuable in resolving several open problems and conjectures in both commutative and noncommutative ring theory. One of the most notable examples is the construction of a chain ring $R$ that is not a factor of a valuation domain (see \cite[Theorem 3.5]{FS1} and \cite[X.6]{FS2}). This ring $R$ is the trivial ring extension of a valuation domain $D$ by a non-standard uniserial divisible $D$-module, providing a counterexample to a question posed by Kaplansky. For further details, see \cite{AJK, BKM, C, KM}.

Note that $A \ltimes M$ is a graded ring with $(A \ltimes M)_0 = A \ltimes 0$, $(A \ltimes M)_1 = 0 \ltimes M$, and $(A \ltimes M)_n = 0$ for $n \geq 2$. Recall from \cite[Theorem 3.3 (1)]{AW} that the homogeneous ideals of $A \ltimes M$ have the form $I \ltimes N$, where $I$ is an ideal of $A$, $N$ is a submodule of $M$, and $IM \subseteq N$. These ideals play a crucial role in studying the properties of $A \ltimes M$ and establishing connections between the properties of $A$ and $M$. In this context, it is particularly interesting that the first example of an ideal that is not homogeneous appears in \cite[Example 2.5]{KM}.

It has been shown that a principal ideal $(a, m)R$ is homogeneous if and only if
\[
(a, m)R = A a \ltimes (A m + aM) = (a, 0)R + (0, m)R,
\]
and every ideal of $A \ltimes M$ is homogeneous if and only if every principal ideal of $A \ltimes M$ is homogeneous \cite[Theorem 3.3 (1)]{AW}. For an integral domain $A$ and an $A$-module $M$, every ideal of $A \ltimes M$ is homogeneous if and only if $M$ is divisible. In this case, all ideals of $A \ltimes M$ have the form $I \ltimes M$ or $0 \ltimes N$, where $I$ is an ideal of $A$ and $N$ is a submodule of $M$ (see \cite[Corollary 3.4 (1)]{AW}). Furthermore, prime (respectively, maximal) ideals of $A \ltimes M$ take the form $P \ltimes M$, where $P$ is a prime (respectively, maximal) ideal of $A$ \cite[Theorem 3.2]{AW}.

In this paper, we further investigate $S$-prime ideals in the trivial ring extension. We begin by providing three counterexamples. Specifically:
\begin{enumerate}
    \item We first show that $S$-prime (resp., $S$-maximal) ideals in the trivial ring extension do not necessarily have the form $P \ltimes M$, where $P$ is an $S_0$-prime (resp., $S_0$-maximal) ideal of $A$ (see Example \ref{ex1}).
    \item In Example \ref{ex2}, we construct an $S$-maximal (and hence $S$-prime) ideal in the trivial ring extension that is not homogeneous.
    \item We provide a non-$S$-prime ideal $J$ in the trivial ring extension $A \ltimes M$ whose projection onto $A$ is an $S_0$-maximal (and hence $S_0$-prime) ideal of $A$.
\end{enumerate}

Moreover, we characterize conditions under which an ideal $J$ (not necessarily homogeneous) in the trivial ring extension is $S$-prime (resp., $S$-maximal). This leads to the following characterization: An ideal $J$ of $R = A \ltimes M$ is $S$-prime if and only if $J_0$ is $S_0$-prime and $M/J_1$ is uniformly $S_0$-torsion, where
\[
J_0 = \{ x \in A \mid (x,0) \in J \}, \quad J_1 = \{ m \in M \mid (0,m) \in J \}
\]
(see Theorems \ref{th1} and \ref{th2}).

As a corollary, we establish that all $S$-prime ideals of $R = A \ltimes M$ take the form $P \ltimes M$, where $P$ is an $S_0$-prime ideal of $A$, if and only if $M$ is $S_0$-divisible (see Corollaries \ref{imp} and \ref{rec}). 

As an application, we explore the transfer of the concepts of compactly $S$-packed rings, coprimely $S$-packed rings and   $S$-$pm$-rings to the trivial ring extension (see theorems \ref{th3}, \ref{th4} and \ref{th5}.

\section{$S$-Prime and $S$-maximal ideals in trivial ring extensions}

Recall from \cite[Theorem 3.8]{AW} that there is a one-to-one correspondence between the saturated multiplicative subsets of $A$ and those of $A \ltimes M$, given by $S \leftrightarrow S \ltimes M$. Furthermore, if $S$ is a multiplicative subset of $A$ and $N$ is a submodule of $M$, then $S \ltimes N$ is a multiplicative subset of $R \ltimes M$ with saturation $(S \ltimes N)^* = S^* \ltimes M$. Moreover, if $S$ is a multiplicative subset of $R = A \ltimes M$ and $S_0$ is its projection onto $A$, then the multiplicative subsets $S$, $S_0 \ltimes N$ (where $N$ is a submodule of $M$), and $S_0 \ltimes M$ all have the same saturation, namely $S_0^* \ltimes M$. This justifies why, in this context, authors focus only on multiplicative subsets of the form $S_0 \ltimes M$, where $S_0$ is a multiplicative subset of $A$.

Let $I$ be an ideal of $R$ and let $J$ be a nonempty subset of $R$. The residual of $I$ by $J$ is defined as
\[
(I : J) = \{x \in R \mid xJ \subseteq I\}.
\]
For a singleton $J = \{a\} \subseteq R$, we use the notation $(I : a)$ instead of $(I : \{a\})$.

Recall from \cite[Proposition 1]{HM} that for a commutative ring $R$ with a multiplicative set $S \subseteq R$, an ideal $\mathfrak{p}$ of $R$ that is disjoint from $S$ is $S$-prime if and only if there exists some $s \in S$ such that $(\mathfrak{p} : s)$ is a prime ideal of $R$.

\textbf{Throughout this paper, the element $s$ associated with an $S$-prime ideal $\mathfrak{p}$ will be denoted by $s_{\mathfrak{p}}$.}

To simplify the presentation of examples, we will focus on cases where the trivial ring extension is of the form $\mathbb{Z} \ltimes M$, where $M$ is a $\mathbb{Z}$-module. In this case, the notions of $S$-prime and $S$-maximal coincide, as established by the following lemma.

\begin{lem} \label{pm}
    Let $M$ be a $\mathbb{Z}$-module, and let $R = \mathbb{Z} \ltimes M$ be the trivial ring extension. If $S$ is a multiplicative subset of $R$, then an ideal $P \nsubseteq Nil(R)$ of $R$ is $S$-prime if and only if it is $S$-maximal.
\end{lem}

\proof
Assume that $P$ is an $S$-prime ideal of $R$ such that $P \nsubseteq Nil(R)$. By \cite[Proposition~1]{HM}, there exists some $s \in S$ such that $(P : s)$ is a prime ideal of $R$ with $Nil(R) \subset (P : s)$. Since every prime ideal of $R$ containing $Nil(R)$ is maximal (see \cite[Theorem~3.2(2)]{AW}), it follows that $(P : s)$ is a maximal ideal of $R$. Therefore, $P$ is $S$-maximal by \cite[Theorem~2.3]{KMO}.
\qed \\

Recall that prime (resp., maximal) ideals of $R$ take the form $P \ltimes M$, where $P$ is a prime (resp., maximal) ideal of $A$ by \cite[Theorem 3.2]{AW}. It is worth noting that the authors in \cite{YETY} studied $S$-prime ideals in the trivial ring extension of rings in \cite[Lemma 2]{YETY}. However, in the following example (as well as in Examples \ref{ex3} and \ref{ex4}), we show that this result does not hold in general.

\begin{exam}\label{ex1}{\rm
    Let $R = \mathbb{Z} \ltimes \mathbb{Z}/6\mathbb{Z}$ and let $S = \{ (2,0)^n \mid n \in \mathbb{N} \}$. Then the ideal $J = 0 \ltimes 2\mathbb{Z}/6\mathbb{Z}$ is an $S$-prime ideal of $R$. However, an $S$-prime ideal in the trivial ring extension $R = A \ltimes M$ does not necessarily have the form $P \ltimes M$, where $P$ is an $S_0$-prime ideal of $R$.
    }
\end{exam}

\proof It is clear that $S \cap J = \emptyset$. Since $\mathbb{Z}$ is an integral domain, it follows from \cite[Theorem 3.2]{AW} that $0 \ltimes \mathbb{Z}/6\mathbb{Z}$ is a prime ideal of $R$. Now, let $x, y \in R$ such that $xy \in J \subseteq 0 \ltimes \mathbb{Z}/6\mathbb{Z}$. Thus, we have either $x \in 0 \ltimes \mathbb{Z}/6\mathbb{Z}$ or $y \in 0 \ltimes \mathbb{Z}/6\mathbb{Z}$. Setting $s = (2,0)$, we obtain that $sx \in (2,0)(0 \ltimes \mathbb{Z}/6\mathbb{Z}) = J$ or $sy \in (2,0)(0 \ltimes \mathbb{Z}/6\mathbb{Z}) = J$. Consequently, $J$ is an $S$-prime ideal of $R$. \qed \\

In the following example, we go further by constructing an $S$-maximal (and hence $S$-prime) ideal in the trivial ring extension that is not homogeneous.

\begin{exam}\label{ex2}{\rm
      Let $R = \mathbb{Z} \ltimes \mathbb{Z}/2\mathbb{Z}$ and let $S = \{ (2,0)^n \mid n \in \mathbb{N} \}$. Then the principal ideal $P$ generated by $(6,1)$ is an $S$-maximal ideal of $R$. However, $P$ is not a homogeneous ideal of $R$.
      }
\end{exam}

\proof First, note that
	$$
	\begin{aligned}
		P &= \{ (k,r)(6,1) \mid (k,r) \in R \} \\
		&= \{ (6k,k) \mid k \in \mathbb{Z} \} \\
		&= \{ (12k,0) \mid k \in \mathbb{Z} \} \cup \{ (12k+6,1) \mid k \in \mathbb{Z} \} \quad (*).
	\end{aligned}
	$$
Thus, we have $P \cap S = \emptyset$. Now, take $s = (4,0)$, and we claim that $(P:s) = 3\mathbb{Z} \ltimes \mathbb{Z}/2\mathbb{Z}$.

Let $x = (3a,e) \in 3\mathbb{Z} \ltimes \mathbb{Z}/2\mathbb{Z}$. Since
\[
sx = (4,0)(3a,e) = (12a,0) = (6,1)(2a,0) \in P,
\]
we conclude that $x \in (P:s)$. Thus, $3\mathbb{Z} \ltimes \mathbb{Z}/2\mathbb{Z} \subseteq (P:s)$. Conversely, let $(a,e) \in (P:s)$. Then
\[
(a,e)(4,0) = (4a,0) \in P,
\]
which implies that $a\in 3\mathbb{Z}$  according to $(*)$. Hence, $(a,e) \in 3\mathbb{Z} \ltimes \mathbb{Z}/2\mathbb{Z}$. Therefore, $(P:s) = 3\mathbb{Z} \ltimes \mathbb{Z}/2\mathbb{Z}$ is a prime ideal of $R$ by \cite[Theorem 3.2]{AW}, since $3\mathbb{Z}$ is a prime ideal of $\mathbb{Z}$. Consequently, $P$ is an $S$-prime ideal of $R$ by \cite[Proposition 1]{HM}, and so $P$ is $S$-maximal by Lemma \ref{pm}.

To demonstrate that $P$ is not a homogeneous ideal of $R$, it is sufficient to show that $(0,1) \notin P$, as stated in \cite[Theorem 3.3 (3)]{AW}. This follows directly from $(*)$. \qed \\

Next, we provide an example of a non-$S$-prime ideal $J$ in the trivial ring extension $A \ltimes M$. However, its projection onto $A$ is an $S_0$-maximal (and hence $S_0$-prime) ideal of $A$.

\begin{exam}\label{ex3}{\rm
    Let $R = \mathbb{Z} \ltimes \mathbb{Z}/2\mathbb{Z}$ and let $S = \{ (3,0)^n \mid n \in \mathbb{N} \}$. Then the principal ideal $P$ generated by $(2,1)$ is not an $S$-prime ideal of $R$. However, $\pi_A(P) = 2\mathbb{Z}$ is an $S_0$-prime ideal of $A$.
    }
\end{exam}

\proof As in Example \ref{ex2}, we have
\[
P = \{ (4k,0) \mid k \in \mathbb{Z} \} \bigcup \{ (4k+2,1) \mid k \in \mathbb{Z} \} \quad (*).
\]
Thus, we see that
\[
(2,0)(2,0) = (4,0) = (2,0)(2,1) \in P.
\]
However, for every $s = (3,0)^n \in S$, we have
\[
s(2,0) = (3^n 2,0) \notin P,
\]
according to $(*)$. Therefore, $P$ is not $S$-prime. However, $\pi_A(P) = 2\mathbb{Z}$ is an $S_0$-maximal ideal of $A$. \qed \\

\begin{lem}\label{sat}
 Let $R$ be a commutative ring with identity, $S$ a multiplicative subset of $R$ with saturation $S^*$, and $I$ an ideal of $R$. Then $I$ is an $S$-prime ideal of $R$ if and only if $I$ is an $S^*$-prime ideal of $R$.
\end{lem}

\proof Assume that $I$ is an $S$-prime ideal of $R$. Then $I \cap S = \emptyset$ and there exists $s \in S$ such that for all $a, b \in R$ with $ab \in I$, we have $sa \in I$ or $sb \in I$. To show that $I$ is an $S^*$-prime ideal of $R$, it is enough to verify that $I \cap S^* = \emptyset$ since $s \in S \subseteq S^*$.

Suppose, for contradiction, that there exists $t \in S^* \cap I$. Since $t \in S^*$, there exists $x \in R$ such that $tx \in S$. This implies that $tx \in S \cap I = \emptyset$, which is a contradiction.

Conversely, assume that $I$ is an $S^*$-prime ideal of $R$. Then $I \cap S^* = \emptyset$, and by \cite[Proposition 1]{HM}, there exists $t \in S^*$ such that $(I : t)$ is a prime ideal of $R$. Since $I \cap S \subseteq I \cap S^*$, we conclude that $I \cap S = \emptyset$.

On the other hand, since $t \in S^*$, there exists $r \in S^*$ such that $s = tr \in S$. As $r \notin (I : t)$ and $(I : t)$ is a prime ideal of $R$, it follows that $(I : s) = (I : t)$. Hence, $I$ is disjoint from $S$ and there exists $s \in S$ such that $(I : s)$ is a prime ideal of $R$. Therefore, $I$ is an $S$-prime ideal of $R$ by \cite[Proposition 1]{HM}. \qed 
\begin{lem}\label{S-P}(\cite[Theorem 1.4.7]{WK})
     Let $I$ be an ideal of a ring $R$, and let $S$ be a multiplicative subset of $R$ with $I \cap S = \emptyset$. Then there exists a prime ideal $\mathfrak{p}$ of $R$ such that $I \subseteq \mathfrak{p}$ and $\mathfrak{p} \cap S = \emptyset$.
\end{lem}

\begin{lem}\label{smax}
    Let $R$ be a  ring, $S$ be a multiplicative subset of $R$ with saturation $S^*$. Let $\p$ be an $S_0$-maximal ideal of $R$, then $t\p$ is an $S$-maximal ideal of $R$ for every $t\in S^*$.
\end{lem}
\proof As $t\p \subseteq \p$ and $\p \cap S= \emptyset$, then $s\p \cap S =\emptyset$. Now let $q$ be an ideal of $R$ disjoint with $S$ such that $t\p \subseteq q $. By Lemma \ref{S-P} there exists a prime ideal $m$ of $R$  disjoint with $S$ such that 
$Q \subseteq m$. Thus  $t\p \subseteq q \subseteq m $. As $m$ is disjoint with $S$, then  $t\notin m$ (if not, there exists $k\in R$ such that $kt\in S\cap m$, a contradiction). Thus $\p\subseteq m$ since $m$ is a prime ideal of $R$ and $t\notin m$. Now as $q$ is an $S$-maximal ideal of $R$ and $m$ is disjoint with $S$. Thus $s_\p  m\subseteq \p $. Take $k\in R$ such that $tk\in S$ and set $s=s_\p tk\in S$. Whence $$sQ\subseteq sm \subseteq tk\p\subseteq t\p.$$ 
Thus $s\p$ is an $S$-maximal ideal of $R$.\qed\\

Our next aim is to characterize when an ideal $J$ of $R = A \ltimes M$ is $S$-prime. To achieve this, we introduce the following notations, which are motivated by the grading of $R$.

Let $A$ be a ring and $M$ an $A$-module. Define the canonical projection maps:
\[
\pi_A: A \ltimes M \to A, \quad (a, x) \mapsto a, \quad \text{and} \quad \pi_M: A \ltimes M \to M, \quad (a, x) \mapsto x.
\]
For an ideal $J$ of $R$, we define:
\[
J_0 = \pi_A(J \cap R_0) = \{ a \in A \mid (a,0) \in J \}
\]
and
\[
J_1 = \pi_M(J \cap R_1) = \{ x \in M \mid (0,x) \in J \}.
\]
It is easy to see that $J_0$ is an ideal of $A$ and $J_1$ is a submodule of $M$.

Note that if $a \in J_0$ and $x \in M$, then $(a,0) \in J$ and $(0,x) \in R$. Consequently,
\[
(a,0)(0,x) = (0, ax) \in J,
\]
which implies that $ax \in J_1$. Hence, $J_0 M \subseteq J_1$, and thus $J_0 \ltimes J_1$ is a homogeneous ideal of $R$ contained in $J$.

It follows that $J$ is a homogeneous ideal of $R$ if and only if $J = J_0 \ltimes J_1$.

Recall from \cite{Z} that an $R$-module $M$ is called a {\em uniformly $S$-torsion module} if there exists an element $s \in S$ such that $sM = 0$.

In the following theorems, we determine the necessary and sufficient conditions for an ideal $J$ (not necessarily homogeneous) in the trivial ring extension to be $S$-prime (resp., $S$-maximal).

\begin{thm}\label{th1}
    Let $A$ be a ring and $M$ an $A$-module. Set $R = A \ltimes M$ as the trivial ring extension, and let $S$ be a multiplicative subset of $R$. Then an ideal $J$ of $R$ is $S$-prime if and only if $J_0$ is $S_0$-prime and $M/J_1$ is uniformly $S_0$-torsion.
\end{thm}

\proof Assume that $J$ is an $S$-prime ideal of $R$. Since the multiplicative subsets $S$ and $S_0 \ltimes 0$ have the same saturation ($=S_0^* \ltimes M$), it follows from Lemma \ref{sat} that $J$ is $(S_0 \ltimes 0)$-prime in $R$. Define $s_J = (s_0,0)$ where $s_0 \in S_0$. Since $J$ is disjoint from $S_0 \ltimes 0$, it follows that $J_0$ is disjoint from $S_0$.

Now, let $a, b \in A$ such that $ab \in J_0$. Then $(ab,0) \in J$, which implies that
\[
(ab,0) = (a,0)(b,0) \in J.
\]
Thus, either
\[
(s_0,0)(a,0) = (s_0a,0) \in J \quad \text{or} \quad (s_0,0)(b,0) = (s_0b,0) \in J.
\]
Therefore, $s_0 a \in J_0$ or $s_0 b \in J_0$, meaning that $J_0$ is an $S_0$-prime ideal of $A$.

On the other hand, for every $x \in M$, we have
\[
(0,x)(0,x) = (0,0) \in J,
\]
which implies that
\[
(s_0,0)(0,x) = (0,s_0x) \in J,
\]
i.e., $s_0 x \in J_1$. Thus, $s M \subseteq J_1$, meaning that $M/J_1$ is uniformly $S_0$-torsion (with respect to $s_0$).

Conversely, assume that $J_0$ is $S_0$-prime and that $M/J_1$ is uniformly $S_0$-torsion. Let $s_2 \in S$ such that $s_2 M \subseteq J_1$.

We first claim that $J \cap S = \emptyset$. Suppose otherwise, i.e., there exists $(s_0,t) \in J \cap S$. Then $s_0 \in S_0$, and we compute
\[
(s_2,0)(s_0,t) = (s_0s_2, s_2t) \in J.
\]
Since $s_2 t \in s_2 M \subseteq J_1$, it follows that $(0, s_2 t) \in J$. Consequently,
\[
(s_0s_2,0) = (s_2,0)(s_0,t) - (0,s_2t) \in J,
\]
which implies that $s_0 s_2 \in J_0$. This contradicts the assumption that $s_0 s_2 \in J_0 \cap S_0 = \emptyset$.

Now, set $s_1 = s_{J_0}$ and define $s = s_1 s_2$. Let $(a,x), (b,y) \in R$ such that
\[
(a,x)(b,y) = (ab, ay+bx) \in J.
\]
Then,
\[
(s_2,0)(a,x)(b,y) = (s_2,0)(ab,ay+bx) = (sab, s_2(ay+bx)) \in J.
\]
Since $s_2(ay+bx) \in sM \subseteq J_1$, we obtain $(0,s_2(ay+bx)) \in J$. Consequently,
\[
(s_2ab,0) = (s_2,0)(a,x)(b,y) - (0,s_2(ay+bx)) \in J,
\]
which implies that $s_2 ab \in J_0$.

Since $J_0$ is $S_0$-prime with respect to $s_1$, and $s = s_1 s_2 \notin J_0$, we conclude that either $s_1 a \in J_0$ or $s_1 b \in J_0$. This means that either
\[
sa \in J_0 \quad \text{or} \quad sb \in J_0,
\]
which implies that either $(sa,0) \in J$ or $(sb,0) \in J$.

On the other hand, since $sM \subseteq s_2M \subseteq J_1$, we also have $sx \in J_1$ and $sy \in J_1$, which leads to $(0,sx) \in J$ and $(0,sy) \in J$. Therefore, we conclude that either
\[
(s,0)(a,x) = (sa,sx) = (sa,0) + (0,sx) \in J
\]
or
\[
(s,0)(b,y) = (sb,sy) = (sb,0) + (0,sy) \in J.
\]
Thus, $J$ is $(S_0 \ltimes 0)$-prime, and so it is $S$-prime by Lemma \ref{sat}. \qed

\begin{thm}\label{th2}
    Let $A$ be a ring and $M$ an $A$-module. Set $R = A \ltimes M$ as the trivial ring extension, and let $S$ be a multiplicative subset of $R$. Then an ideal $J$ of $R$ is $S$-maximal if and only if $J_0$ is $S_0$-maximal and $M/J_1$ is uniformly $S_0$-torsion.
\end{thm}

\proof Let $J$ be an $S$-maximal ideal of $R$. Then $J$ is also $S$-prime by \cite[Proposition 10]{szariski}. Thus, $J_0$ is $S_0$-prime and $M/J_1$ is uniformly $S_0$-torsion. It remains to show that $J_0$ is an $S_0$-maximal ideal of $A$.

Let $I$ be an ideal of $A$ disjoint from $S_0$ such that $J_0 \subseteq I$. Then we can easily see that $(s_1,0)J \subseteq I \ltimes M$ for some $s_1 \in S$ such that $s_1M \subseteq J_1$. Since $S \cap (I \ltimes M) = \emptyset$ and $(s_1,0)J$ is also $S$-maximal by Lemma \ref{smax}, it follows that
\[
(s_0,t)(I \ltimes M) = s_0 I \ltimes (s_0 M + It) \subseteq (s_0,0)J \subseteq J
\]
for some $(s_0,t) \in S$. Thus, $s_0 I \subseteq J_0$, implying that $J_0$ is an $S_0$-maximal ideal of $A$.

Conversely, assume that $J_0$ is $S_0$-maximal and that $M/J_1$ is uniformly $S_0$-torsion. By \cite[Theorem 2.3]{KMO}, there exists $s_1 \in S_0$ such that $(J_1:s_1)$ is a maximal element in the set of all ideals of $A$ disjoint from $S_0$. Let $s_2 \in S_0$ be such that $s_2M \subseteq J_1$, and define $s = s_1s_2$. Since
\[
(J_1:s_1) \subseteq (J_1:s) \subseteq (J_1:s^2),
\]
and since $(J_1:s^2)$ is an ideal of $A$ disjoint from $S_0$, it follows that $(J_1:s_1) = (J_1:s) = (J_1:s^2)$ is a maximal element of the set of all ideals of $A$ disjoint from $S_0$.

Now, let $t \in M$ such that $(s,t) \in S$. We claim that
\[
(J:(s,t)^2) = (J_0:s) \ltimes M.
\]
First, note that for every $x \in M$, since $sx \in J_1$, we have
\[
(0,x)(s,t) = (0,x) \in J,
\]
which implies that $0 \ltimes M \subseteq (J:(s,t)^2)$. Also, for every $a \in (J_0:s)$, since $sa \in J_0$, it follows that
\[
(a,0)(s,t)(s,t) = (as,at)(s,t) = (as^2, 2ast) = (as^2,0) + (0,s(2at)) \subseteq J.
\]
Since $s(2at) \in sM \subseteq J_1$ and $as^2 \in J_0$, we obtain
\[
(J_0:s) \ltimes 0 \subseteq (J:(s,t)^2).
\]
Thus,
\[
(J_0:s) \ltimes M = (J_0:s) \ltimes 0 + 0 \ltimes M \subseteq (J:(s,t)^2).
\]

For the reverse inclusion, let $(a,x) \in (J:(s,t)^2)$. Then
\[
(a,x)(s,t)^2 = (a,x)(s^2,2st) = (as^2, 2ast+s^2x) \in J.
\]
Thus,
\[
(as^2,0) = (as^2, 2ast+s^2x) - (0,s(2at+sx)) \in J
\]
since $s(2at+sx) \in J_1$. Hence, $as^2 \in J_0$, and so $x \in (J_0:s^2) = (J_0:s)$. Therefore,
\[
(J:(s,t)^2) = (J_0:s) \ltimes M.
\]

To show that $J$ is an $S$-maximal ideal of $R$, it suffices to prove that $(J:(s,t)^2) = (J_0:s) \ltimes M$ is a maximal element in the set of all ideals of $R$ disjoint from $S$. For this, let $K$ be an ideal of $R$ disjoint from $S$ such that
\[
(J:(s,t)^2) = (J_0:s) \ltimes M \subseteq K.
\]
Since $0 \ltimes M \subseteq (J_0:s) \ltimes M \subseteq K$, it follows that $K = I \ltimes M$ for some ideal $I$ of $A$ by \cite[Theorem 3.1]{AW}. On the other hand, since $K = I \ltimes M$ is disjoint from $S$, we have that $I$ is disjoint from $S_0$. By the maximality of $(J_0:s)$, it follows that $(J_0:s) = I$ (noting that $(J_0:s) \subseteq I$ follows by projecting onto $A$ via the inclusion $(J_0:s) \ltimes M \subseteq K = I \ltimes M$). Therefore,
\[
(J:(s,t)^2) = (J_0:s) \ltimes M = I \ltimes M = K.
\]
\qed

\begin{cor}
     Let $A$ be a ring and $M$ an $A$-module. Set $R = A \ltimes M$ as the trivial ring extension, and let $S$ be a multiplicative subset of $R$. Let $J$ be an ideal of $R$. Then the following assertions are equivalent:
     \begin{enumerate}
         \item $J$ is an $S$-prime (resp., $S$-maximal) ideal of $R$;
         \item $J_0 \ltimes J_1$ is an $S$-prime (resp., $S$-maximal) ideal of $R$;
         \item $J_0$ is $S_0$-prime (resp., $S_0$-maximal) and $M/J_1$ is uniformly $S_0$-torsion.
     \end{enumerate}
\end{cor}

\begin{cor}
     Let $A$ be a ring and $M$ an $A$-module. Set $R = A \ltimes M$ as the trivial ring extension, and let $S$ be a multiplicative subset of $R$. Let $I$ be an ideal of $A$ and $N$ a submodule of $M$ such that $IM \subseteq N$. Then the following assertions are equivalent:
     \begin{enumerate}
         \item $I \ltimes N$ is an $S$-prime (resp., $S$-maximal) ideal of $R$;
         \item $I$ is $S_0$-prime (resp., $S_0$-maximal) and $M/N$ is uniformly $S_0$-torsion.
     \end{enumerate}
\end{cor}

Let $S$ be a multiplicative subset of a ring $R$, and let $M$ be an $R$-module. Recall from \cite[Definition 1.6.10]{WK} that $M$ is said to be {\em $S$-divisible} if for any $s \in S$ and any $x \in M$, there exists $y \in M$ such that $sy = x$, equivalently $sM = M$.

\begin{cor}\label{imp}
     Let $A$ be a ring and $M$ an $S_0$-divisible $A$-module. Set $R = A \ltimes M$ as the trivial ring extension, and let $S$ be a multiplicative subset of $R$. Then
     all $S$-prime (resp., $S$-maximal) ideals of $R$ have the form $P \ltimes M$, where $P$ is an $S_0$-prime (resp., $S_0$-maximal) ideal of $A$.
\end{cor}

\proof Let $P$ be an $S$-prime ideal of $R$. By Theorem \ref{th1}, $P_0$ is $S_0$-prime and $M/P_1$ is uniformly $S_0$-torsion. Thus, there exists $s \in S$ such that $sM \subseteq P_1 \subseteq M$. Since $M$ is an $S_0$-divisible $A$-module, it follows that $sM = M$, and hence $P_1 = M$. Therefore, $J = P_0 \ltimes M$, where $P_0$ is an $S_0$-prime ideal of $A$. The same reasoning applies to $S$-maximal ideals. \qed \\

When $S$ is a subset of units of $R$, we recover the following classical result.

\begin{cor}(\cite[Theorem 3.2]{AW})
     Let $A$ be a ring and $M$ an $A$-module. Set $R = A \ltimes M$ as the trivial ring extension. Then all prime (resp., maximal) ideals of $R$ have the form $P \ltimes M$, where $P$ is a prime (resp., maximal) ideal of $A$.
\end{cor}

\begin{cor}\label{rec}
     Let $A$ be a ring and $M$ an $A$-module. Set $R = A \ltimes M$ as the trivial ring extension, and let $S$ be a multiplicative subset of $R$. If
     all $S$-prime ideals of $R$ have the form $P \ltimes M$, where $P$ is an $S_0$-prime ideal of $A$, then $M$ is $S_0$-divisible.
\end{cor}

\proof Since $S_0$ is a multiplicative subset of $A$, Lemma \ref{S-P} guarantees the existence of a prime ideal $P$ of $A$ disjoint from $S_0$. Consequently, $P \ltimes M$ is a prime ideal of $R$ disjoint from $S_0 \ltimes 0$, and hence it is an $(S_0 \ltimes 0)$-prime ideal of $R$.

Let $s \in S$. Since $(s,0) \in S_0 \ltimes 0$, it follows from \cite[Proposition 2 (1)]{HM} that
\[
(s,0)P \ltimes M = sP \ltimes sM
\]
is an $S$-prime ideal of $R$ by combining Lemma \ref{sat} and \cite[Theorem 3.8]{AW}. By hypothesis, this implies that $sM = M$. Therefore, $M$ is $S_0$-divisible. \qed

\begin{rem}{\rm
    With an alternative proof, we can also replace the statement in the previous corollary that
    ``$S$-prime ideals of $R$ have the form $P \ltimes M$, where $P$ is an $S_0$-prime ideal of $A$"
    with the analogous result for $S$-maximal ideals:
   ``$S$-maximal ideals of $R$ have the form $p \ltimes M$, where $P$ is an $S_0$-maximal ideal of $A$."
}
\end{rem}

\begin{cor}\label{S-torsion}
     Let $A$ be a ring and $M$ a uniformly $S_0$-torsion $A$-module. Set $R = A \ltimes M$ as the trivial ring extension, and let $S$ be a multiplicative subset of $R$. Then an ideal $J$ of $R$ is $S$-prime (resp., $S$-maximal) if and only if $J_0$ is an $S_0$-prime (resp., $S_0$-maximal) ideal of $A$.
\end{cor}

We conclude this section with the following example, which exhibits an $S$-prime ideal in the trivial ring extension that is not homogeneous.

\begin{exam}\label{ex4}{\rm
    Let $R = \mathbb{Z} \ltimes \mathbb{Z}/4\mathbb{Z}$ and set $S = \{ (2,0)^n \mid n \in \mathbb{N} \}$. Then the principal ideal $J = (6,1)R$ is an $S$-prime ideal of $R$ that is not homogeneous.
    }
\end{exam}

\proof First, note that $n \in J_0$ if and only if $(n,0) \in J$, which means there exists $(m,x) \in R$ such that
\[
(n,0) = (6,1)(m,x) = (6m,6x+m).
\]
Thus, $J_0 \subseteq 6\mathbb{Z}$. However, observe that $6 \notin J_0$, since if $6 \in J_0$, then there would exist $(m,x) \in R$ such that
\[
(6,0) = (6,1)(m,x) = (6m,6x+m).
\]
This forces $m = 1$, leading to the equality $6x + 1 = 0$, i.e., $6x + 1 \in 4\mathbb{Z}$, which is a contradiction.

On the other hand, since
\[
(12,0) = (6,1)(2,1) \in J,
\]
it follows that $12 \in J_0$. Consequently, we conclude that $J_0 = 12\mathbb{Z}$, which is disjoint from $S_0 = \{ 2^n \mid n \in \mathbb{N} \}$.

Now, the fact that $(J_0:4) = 3\mathbb{Z}$ is a prime ideal of $A$ implies that $J_0$ is an $S_0$-prime ideal by \cite[Proposition 1]{HM}. Hence, $J$ is $S$-prime according to Corollary \ref{S-torsion}.

To show that $J$ is not a homogeneous ideal of $R$, it is sufficient to verify that $(6,0) \notin J$, as stated in \cite[Theorem 3.3 (3)]{AW}. However, this has already been established. \qed

\section{Transfer of  compactly $S$-packed, coprimely $S$-packed and   $S$-$pm$-Property to Trivial Ring Extension}

A ring $R$ is called a compactly packed ring if whenever $Q \subseteq \bigcup_{i \in \Gamma} \p_i$ for some ideal $Q$ of $R$ and a family of prime ideals $\left\{\p_i\right\}_{i \in \Gamma}$ of $R$, then there exists $i \in \Gamma$ such that $Q \subseteq \p_i$ \cite{RV}.  On the other hand, Erdoğdu \cite{E} presented the notion of coprimely packed rings, a generalization of compactly packed rings. A ring $R$ is called a coprimely packed ring if whenever $Q+\p_i=R$ for some ideal $Q$ of $R$ and a family of prime ideals $\left\{\p_i\right\}_{i \in \Gamma}$ of $R$, then $Q \nsubseteq \bigcup_{i \in \Gamma} \p_i$. Note that every compactly packed ring is a coprimely packed ring. The converse of this indication is correct when $R$ is a domain with a Krull dimension is one \cite[Proposition 2.2]{E}. In \cite{YETY}, Yavuz et al. introduced the notions of compactly $S$-packed rings and  coprimely $S$-packed rings as follows:
\begin{defn} Let $R$ be a ring and $S$ be a multiplicative subset of $R$. Then
  \begin{enumerate}  \item (\cite[Definition 2]{YETY}) $R$ is called a compactly $S$-packed ring if $Q \subseteq \bigcup_{i \in \Gamma} P_i$, where $Q$ is an ideal of $R$ and $\left\{P_i\right\}_{i \in \Gamma}$ is a family of S-prime ideals of $R$, then there exists $k \in \Gamma$ and $s \in S$ such that $s Q \subseteq P_k$.
     \item  (\cite[Definition 3]{YETY}) $R$ is called a coprimely $S$-packed ring if whenever $Q+P_i = R$ for some ideal $Q$ of $R$ and a family of $S$-prime ideals $\left\{P_i\right\}_{i \in \Gamma}$ of $R$, then $s Q \nsubseteq \bigcup_{i \in \Gamma} P_i$ for every $s \in S$.
\end{enumerate}
\end{defn}

Among other things,  they provided the transfer of these notions  to the trivial extension of the ring  (cf. \cite[Theorems 8 and 10]{YETY}). However, their proofs rely on  \cite[Lemma 2]{YETY}, which states that all $ S $-prime ideals of $ A\ltimes  M $ have the form $ P\ltimes  M $ where $P$ is an $S_0$-prime ideal of $A$. This statement is not generally holds (see Examples \ref{ex1}, \ref{ex2} and \ref{ex4}). Nevertheless, in the specific case where $ M $ is $S_0$-divisible, \cite[Lemma 2]{YETY} holds (see Corollary \ref{imp}), and as a result,  \cite[Theorems 8 and 10]{YETY}  are also valid in this context. 
This raises a natural question: do these theorems remain valid in the general case? The aim of this section is to investigate the conditions under which the trivial extension is compactly $S$-packed ring (resp., coprimely $S$-packed ring). As a corollary, we will demonstrate the validity of \cite[Theorems 8 and 10]{YETY}.

\begin{thm}\label{th3}
Let  $M$ be an $A$-module and let $S$ be a multiplicative subset of the trivial ring extension   $R= A\ltimes  M$. Then the following statements are equivalent:
\begin{enumerate}
    \item 
 $ A$ is a compactly $S_0$-packed ring;
  \item  $R$ is a compactly $S$-packed ring.
\end{enumerate}
\end{thm}
\proof   $1)  \Rightarrow 2)$ Assume that  $A$ is a compactly $S_0$-packed ring. Suppose that $Q \subseteq \bigcup_{i \in \Gamma} P_i$ for some ideal $Q$ of $A\ltimes  M$ and a family of $S$-prime ideals  $\left\{P_i\right\}_{i \in \Gamma}$ of $A\ltimes   M$. Take $\p_i=(P_i:s_{P_i})$ for every $i \in \Gamma$. Then $\left\{\p_i\right\}_{i \in \Gamma}$ is family of prime ideals
 of $R$ disjoint with $S$. Whence for every $i \in \Gamma $, $\p_i=I_i\ltimes  M$ for some prime ideal $I_i$ of $A$ by \cite[Theorem 3.2]{AW}. Now the fact that $\p_i=I_i\ltimes  M$ is disjoint with $S$ for all $i\in \Gamma$, implies that $I_i\cap S_0=\emptyset$. Thus all $I_i$ are prime ideals of $A$ disjoint with $S_0$ and so all  $I_i$  are $S_0$-prime ideal of $A$. On the other hand we have  $$Q \subseteq \bigcup_{i \in \Gamma} P_i \subseteq \bigcup_{i \in \Gamma} \p_i =\bigcup_{i \in \Gamma} I_i \ltimes  M= (\bigcup_{i \in \Gamma} I_i)\ltimes  M.$$ 
Thus, by projection into  $A$ we conclude that $$q:=\pi_A(Q)\subseteq \bigcup_{i \in \Gamma} I_i =\pi_A(\bigcup_{i \in \Gamma} \p_i).$$   Since $A$ is a compactly $S$-packed ring, there exists $s \in S$ such that $s q \subseteq I_i$ for some $i \in \Gamma$. Thus  $$(s, t) q\ltimes  M =sq\ltimes  ( qt+sM) \subseteq sq\ltimes  M  \subseteq I_i \ltimes  M=\p_i.$$
  Note that as $q:=\pi_A(Q)$, then $Q \subseteq q\ltimes  M$. Thus 
  $$(s, t)Q \subseteq (s, t) q\ltimes  M  \subseteq \p_i.$$
  Since $\p_i$ is a prime ideal of $R$ disjoint with $S$ and $(s, t)\in S$ (and so $(s, t)\notin S$), it follows that 
  $Q \subseteq \p_i=(P_i:s_{P_i})$. Thus $$s_{P_i} Q \subseteq s_{P_i} (P_i:s_{P_i})\subseteq P_i.$$
Therefore, $ A \ltimes  M$ is a compactly $S$-packed ring.\par
$2)   \Rightarrow 1)$ Assume that  $A\ltimes  M $ be a compactly $S$-packed ring. Suppose that $q \subseteq \bigcup_{i \in \Gamma} p_i$ for some ideal $q$ of $A$ and a family of $S_0$-prime ideals $\left\{p_i\right\}_{i \in \Gamma}$ of $A$. Then, $q\ltimes  M \subseteq \bigcup_{i \in \Gamma}\left(p_i\ltimes  M\right)$ where $\left(p_i \ltimes  M\right)$ is $S$-prime ideal  of $A\ltimes  M$ for all $i \in \Gamma$ according to Theorem \ref{th1}. Since $A\ltimes  M$ is a compactly $S$-packed ring, there exists $(s, t) \in S$ such that
$$sq\ltimes (sM+qm) =(s, m)(q\ltimes  M)  \subseteq p_i\ltimes  M$$
for some $i\in \Gamma$. Thus, by projection into $A$ we conclude that 
$sq\subseteq p_i$ and $s\in S_0$. Thus $ A$ is a compactly $S_0$-packed ring.\qed
\begin{cor}(\cite[Theorem 8]{YETY})
     Let  $M$ be  an $A$-module and $S$ be a multiplicative subset of  $A$. Then the following statements are equivalent:
\begin{enumerate}
    \item  $A$ is a compactly $S_0$-packed ring;
\item  $ R$ is a compactly $(S\ltimes  0)$-packed ring;
\item  $ R$  is a compactly $(S\ltimes  M)$-packed ring.
\end{enumerate}
\end{cor}
When $S$ is a subset of units of $R$,  we obtain the following result:

\begin{cor}
    Let  $M$ be an $A$-module. Then the following statements are equivalent:
\begin{enumerate}
    \item  $A$ is a compactly packed ring;
\item  $ R$  is a compactly packed ring.
\end{enumerate}
\end{cor}

\begin{thm}\label{th4}
Let $M$ be an $A$-module and $S$ be a  multiplicative subset of  the trivial ring extension $R= A\ltimes  M$. Then the following statements are equivalent:
\begin{enumerate}
    \item 
 $ A$ is a coprimely $S_0$-packed ring;
  \item  $R$ is a coprimely $S$-packed ring.
\end{enumerate}
\end{thm}
\proof   $1)  \Rightarrow 2)$  Assume that  $R$  is a coprimely $S_0$-packed ring. Suppose that $Q+P_i=A\ltimes   M$ for some ideal $Q$ of $A\ltimes  M$ and a family of $S$-prime ideals $\left\{P_i\right\}_{i \in \Gamma}$ of $A\ltimes  M$. Take $\p_i=(P_i:s_{P_i})$ for every $i \in \Gamma$. Then $\left\{\p_i\right\}_{i \in \Gamma}$ is family of prime ideals
 of $R$ disjoint with $S$. Whence for every $i \in \Gamma $, $\p_i=I_i\ltimes  M$ for some prime ideal $I_i$ of $A$ by \cite[Theorem 3.2]{AW}. Now the fact that $\p_i=I_i\ltimes  M$ is disjoint with $S$ for all $i\in \Gamma$, implies that $I_i\cap S_0=\emptyset$. Thus all $I_i$ are prime ideals of $A$ disjoint with $S_0$ and so all  $I_i$  are $S_0$-prime ideal of $R$.
   Take $J=\pi_A(Q)$, then $Q\subseteq J \ltimes  M$. Therefore $$A\ltimes  M= Q+P_i\subseteq Q+ \p_i \subseteq J \ltimes  M + I_i\ltimes  M= (J+I_i) \ltimes  M$$ for all $i \in \Gamma$. 
   Thus, $J+I_i=A$. Because $A$ is a coprimely $S_0$-packed ring and   all  $I_i$  are $S_0$-prime ideal of $A$, we obtain that $s_0 J \nsubseteq \bigcup_{i \in \Gamma} I_i$ for all $s_0 \in S_0$. Now, we claim that  $(s, t) Q \nsubseteq \bigcup_{i \in \Gamma} P_i$ for all $(s, m) \in S$. If this is not the case, then there exists  $(s, t)\in S$ such that $$(s, t) Q \subseteq \bigcup_{i \in \Gamma} P_i \subseteq \bigcup_{i \in \Gamma} (P_i:s_{P_i}) =   \bigcup_{i \in \Gamma}\p_i = \bigcup_{i \in \Gamma} I_i\ltimes  M.$$
   Whence for every $x\in J$, from the description of $J$, there exists $m\in M$ such that $(x,m)\in Q$
  and so $$(s,t)(x,m)=(sx,sm+xt)\in  (s, t) Q  \subseteq \bigcup_{i \in \Gamma} I_i\ltimes  M.$$
  Thus $sx\in \bigcup_{i \in \Gamma} I_i$, that is $sJ\subseteq \bigcup_{i \in \Gamma} I_i$, a contraduction.
Thus, $A\ltimes  M$ is a coprimely $S$-packed ring.\par
$2)   \Rightarrow 1)$ Assume that $A\ltimes   M$ is a coprimely $S$-packed ring. Suppose that $q+p_i=A$ for some ideal $q$ of $A$ and a family of $S_0$-prime ideals $\left\{p_i\right\}_{i \in \Gamma}$ of $A$. Take
$Q=q\ltimes  M$ and $P_i= p_i\ltimes  M$. Then,  $Q+P_i=A\ltimes   M$ and $\left\{P_i\right\}_{i \in \Gamma}$ is a  family of $S$-prime ideals of $A\ltimes  M$ by  Theorem \ref{th1}. Thus, by hypothesis for every $s=(s_0,t)\in S$
$$sq\ltimes  (sM+qM)= (s, t) Q \nsubseteq \bigcup_{i \in \Gamma} P_i= (\bigcup_{i \in \Gamma} p_i)\ltimes  M.$$
Thus, for every $s\in S_0$, $sq \nsubseteq \bigcup_{i \in \Gamma} p_i$. Therefore, $A$ is a coprimely $S_0$-packed ring.\qed
\begin{cor}(\cite[Theorem 10]{YETY})
  Let $M$ be an $A$-module and $S$ is a multiplicative subset of $A$. Then the following statements are equivalent:
\begin{enumerate}
    \item $A$ is a coprimely $S$-packed ring;
\item  $ A\ltimes  M$ is a coprimely $(S\ltimes  0)$-packed ring;
\item  $ A\ltimes  M$  is a coprimely $(S\ltimes  M)$-packed ring.
\end{enumerate}
\end{cor}
When $S$ is a subset of units of $R$,  we obtain the following result
\begin{cor}
Let $M$ be an $A$-module. Then the following statements are equivalent:
\begin{enumerate}
    \item $A$ is a coprimely $S$-packed ring;
\item  $ A\ltimes  M$ is a coprimely packed ring.
\end{enumerate}
\end{cor}

 In \cite{DO}, De Marco and Orsatti introduced the notion of a $pm$-ring as a  ring in which every prime ideal is contained in a unique  maximal ideal. This class of rings contains the class of regular rings, local rings,
 zero-dimensional rings, rings of functions, etc. They studied some properties of the
 prime spectrum $Spec(R)$ and the maximal spectrum $Max(R)$ of a $pm$-ring; and
 various topological and algebraic properties of $pm$-rings have been given  (see, for instance, \cite{DO, MMM}).

Similarly, a ring $ R $ with a multiplicative subset $ S $ is called an $ S $-$pm$-ring if every $ S $-prime ideal of $ R $ is contained in a unique $S$-maximal ideal of $ R $.
\begin{thm}\label{pm}\label{th5}
    Let $M$ be an $A$-module and $S$ be a  multiplicative subset of  the trivial ring extension $R= A\ltimes  M$. Then the following statements are equivalent:
\begin{enumerate}
    \item 
 $ R$ is an  $ S $-$pm$-ring;
  \item  $A$ is an  $ S_0$-$pm$-ring and $M$ is a $S_0$-divisible $A$-module.
\end{enumerate}
\end{thm}
$1)  \Rightarrow 2)$  Assume that $R= A\ltimes  M$ is an  $ S $-$pm$-ring. Let $s\in S_0$ and  $\p$ be an $S_0$-maximal ideal of $A$, then $s\p$ is also an $S_0$-maximal ideal of $R$ by Lemma \ref{smax}. Thus $$(s,0) p\ltimes  M= s\p \ltimes  sM \subseteq \p \ltimes  M.$$ As $s\p \ltimes  sM$ and $\p \ltimes  M$ are $S$-maximal ideal of $R$ according to Theorem \ref{th2} and $ R$ is an  $ S $-$pm$-ring, we conclude that $s\p \ltimes  sM=\p \ltimes  M$. Thus $sM=M$ for every $s\in S_0$, whence $M$ is $S_0$-divisible.  Let $p$ be an $S_0$-prime ideal of $A$ and let $m$ be an $S$-maximal ideal of $A$ such that $p\subseteq m$. Then $p\ltimes  M$  is an $S$-prime ideal of $R$ by Theorem \ref{th1} and $m\ltimes  M$  is an $S$-maximal ideal of $R$ by Theorem \ref{th2} with $p\ltimes  M \subseteq m\ltimes  M$. As $R$ is  $ S $-$pm$-ring, then $m\ltimes  M$ is the only $S$-maximal ideal of $R$   contains the ideal prime ideal $p \ltimes  M$. Thus $m$ is the only $S_0$-maximal ideal of $A$  contains the ideal prime ideal $p$. Thus $A$ is an  $ S $-$pm$-ring.\par
$2) \Rightarrow 1)$ As $M$ is a $S_0$-divisible $A$-module. Then  $$\operatorname{Spec}_S(R)=\{P \ltimes M \mid P \in \operatorname{Spec}_{S_0}(A)\}$$ and  $$\operatorname{Max}_S(R)=\{m \ltimes M \mid m \in \operatorname{Max}_{S_0}(A)\}$$ 
by Corollary \ref{imp}. Thus, it is easy to see that $ R$ is an  $ S $-$pm$-ring since $ A$ is an  $ S _0$-$pm$-ring.\qed\\

It is straightforward to observe that if a ring $ R$ is an $ S$-$ pm $-ring, then $R_S$ is a $ pm$-ring. However, the following example proves that the converse does not hold in general.

\begin{exam}
    Let $R(\Z \ltimes \Z)_{2\Z \ltimes \Z}=\Z_2\ltimes  \Z_2$ and $S=\{ (2,0)^n\mid n\in \N \}$, Then $R$ is a a  $pm$-ring (and so $R_S$) since $R$ is a  local ring. However, $R$ is not an $S$-$pm$-ring acoording to Theorem \ref{pm} since $\Z_2$ is not $S_0$-divisbile where $S_0=\{ 2^n \mid n\in \N \}$.
\end{exam}

\bigskip

\noindent {\bf Funding information:} H. Kim was supported by Basic Science Research Program through the
National Research Foundation of Korea(NRF) funded by the Ministry of Education
(2021R1I1A3047469).

\bigskip

\noindent {\bf Conflict of interest:} The authors state no conflict of interest.

\end{document}